%% file: seriesparallel.tex
\documentclass[letterpaper,12pt,reqno]{amsart}
\usepackage[top=1in,bottom=1in]{geometry}
\usepackage{times}
\usepackage{amsmath}
\usepackage{amssymb}
\usepackage{amsfonts}
\usepackage{amscd}
\usepackage{amsthm}
\usepackage{amsxtra}
\usepackage{stmaryrd}

\usepackage[all]{xy}  \CompileMatrices
\usepackage{mathrsfs}
\usepackage{nicefrac}
\usepackage{bbold}
\usepackage{graphicx}
\usepackage{color}
\usepackage[raggedright,IT,hang]{subfigure}
\usepackage{picins}
\usepackage{algorithmic}

\input{defs.amsart}

\begin{document}

\title{The Cops \& Robber game on series-parallel graphs}
\author{Dirk Oliver Theis}
\address{Dirk Oliver Theis \\
  Service de G\'eometrie Combinatoire et Th\'eorie des Groupes \\
  D\'e\-partement de Math\'ematique \\
  Universit\'e Libre de Bruxelles \\
  Brussels, Belgium}%
\email{theis@uni-heidelberg.de}%
\thanks{Research supported by \textit{Communaut\'e fran\c caise de Belgique -- Actions de Recherche Concert\'ees.}
  Author currently supported by \textit{Fonds National de la Recherche Scientifique, FNRS, Belgium.}}

\date{Thu Jan 17 21:54:20 CET 2008}


\begin{abstract}
  The Cops and Robber game is played on undirected finite graphs.  $k$ cops and one robber are positioned on vertices and take turn in moving along edges.  The
  cops win if, after a move, a cop and the robber are on the same vertex.  A graph is called $k$-copwin, if the cops have a winning strategy.
  It is known that planar graphs are 3-copwin (Aigner \& Fromme, 1984) and that outerplanar graphs are 2-copwin (Clarke, 2002).  In this short note, we prove that
  series-parallel (i.e., graphs with no $K_4$ minor) graphs are 2-copwin.

  It is a well-known trick in the literature of cops \& robber games to define variants of the game which impose restrictions on the possible strategies of the
  cops (see Clarke, 2002).  For our proof, we define the ``cops \& robber game with exits''.

  Our proof yields a winning strategy for the cops.
\end{abstract}
\maketitle

\section{Introduction}

Let $G$ be a finite undirected graph.  For an integer $k\ge 1$, the \textit {$k$-cops and robber game} on $G$ is played as follows: $k$ cops place themselves on
(not necessarily distinct) vertices of the graph; then one robber chooses his vertex.  Now, starting with the $k$ cops, the cops and the robber move alternately,
where a move consists either in passing or in sliding along an edge.  All players have complete information.  The game ends when a cop and the robber are
positioned on the same vertex, i.e., the cops catch the robber, in which case the cops win.  It may proceed forever, though, in which case the robber wins.
We say that $G$ is \textit{$k$-copwin,} if, in the $k$-cops game, the cops have a winning strategy.  The smallest integer $k$ such that a graph is $k$-copwin is
called its \textit{cop-number.}

The class of graphs $1$-copwin graphs has been characterized \cite{NowakowskiWinkler83,QuilliotTh}.  But for every other value of $k$, it appears to be out of
reach of today's tools to give a description of $k$-copwin graphs.  An algorithm recognizing $k$-copwin graphs has been proposed \cite{HahnMacgillivray06}.  For
the case $k=2$, there have been approaches based on restricting the possible strategies the cops may pursue (e.g.
\cite{ClarkePhD,ClarkeNowakowski05,ClarkeNowakowski}).  Moreover, bounds on the cop-number have been investigated (e.g., \cite{BerarducciIntrigila93}), in
particular upper bounds based on coverings with shortest paths (e.g., \cite{FitzpatrickPhD,Clarke}).

It is known that planar graphs are 3-copwin \cite{AignerFromme84} and that outerplanar graphs are 2-copwin \cite{ClarkePhD}.  Surprisingly, nothing is known about
series-parallel graphs, i.e., graphs which do not contain a $K_4$ minor.  In this short note, we prove that series-parallel graphs are 2-copwin.

As mentioned above, it is a known trick in the literature of cops \& robber games to define variants of the game which impose restrictions on the possible
strategies of the cops.  For our proof, we define the ``cops \& robber game with exits''.  Here, a set of ``exit vertices'' is given.  Informally, through an exit
vertex, the robber may ``leave the graph'', and thus win.  We may visualize this as the presence of special edges, which the cops cannot use, from the exits to a
``save haven''.

More formally, if $G$ is an undirected finite reflexive graph, $X\subset V(G)$ a set of vertices and let $k\ge \#X$, the \textit{$k$-cops and robber game on $G$
  with set of exits $X$} is played as follows.  $\#X$ cops are placed on the exit vertices.  The remaining $k-\#X$ cops are placed by the robber (in an adversary
manner), who also selects his own vertex.  Then, starting with the $k$ cops, the cops and the robber move alternately.  Apart from the selection of the starting
vertices, the cops have control over themselves.  The game ends when either a cop and the robber are on the same vertex, in which case the cops win, or when, after
a move of the $k$ cops, the robber stands on an exit vertex alone, in which case the robber wins.  The robber also wins if the game never ends.

We say that a graph $G$ is \textit{$X$-exit $k$-copwin,} if $X\subset V(G)$ and there is a winning strategy for the cops and robber game on $G$ with set of exits
$X$ and $k$ cops.  Clearly, if $G$ is $X$-exit $k$-copwin, then it is $k$-copwin.  We will prove the following.

\begin{theorem}\label{thm:main}
  Every series-parallel graph is $\{x\}$-exit $2$-copwin for each of its vertices $x$.
\end{theorem}

%

In the next section, we will prove Theorem~\ref{thm:main}.  We direct the reader to the fact that our proof yields a winning strategy for the cops, although we
will not write it down explicitly.  

\section{Proof}

In what follows, all series-parallel graphs are assumed to have at least two vertices.
Some terminology will be needed.
We say that a graph $G$ which can be obtained out of a single edge $uv$ by successively duplicating and subdividing edges is a \textit{path-like} series-parallel
graph, and that the vertices of $G$ corresponding to $u$ and $v$ are a pair of \textit{ends} for $G$.  For easy reference, we note the following.

\begin{remark}\label{rem:ex-coend}
  If $G$ is a 2-connected series parallel graph and $u\in V(G)$ arbitrary, then $G$ is path-like and there exists a vertex $v\in V(G)$ such that $u$, $v$ is a pair
  of ends for $G$.
\end{remark}

We will prove the following lemma below.

\begin{lemma}\label{lem:2-conn}
  Let $G$ be  a 2-connected series parallel graph, and let $u$, $v$ be a pair of ends for $G$.  Then $G$ is $\{u,v\}$-exit 2-copwin.
\end{lemma}

With this, we can prove Theorem~\ref{thm:main} by decomposing into blocks.

\begin{proof}[Proof of Theorem~\ref{thm:main}]
  We prove by induction on the number of blocks of the graph.  If $G$ is an edge, then there is nothing to prove; and if $G$ is 2-connected, Lemma~\ref{lem:2-conn}
  together with Remark~\ref{rem:ex-coend} implies that $G$ is $\{x\}$-exit 2-copwin for every vertex $x$.

  Now let $G$ be a series-parallel graph, and $u$ one of its vertices (the exit vertex).  Among the blocks of $G$ which contain $u$ let $B$ be the one with the
  following property:
  \begin{equation}\label{eq:block-to-rob}\tag{$\ast$}%
    \text{$B$ contains an edge of a path from $u$ to the robber.}
  \end{equation}

  Since one of the two cops sits on $u$, condition~\eqref{eq:block-to-rob} holds until the robber reaches~$u$, which we will make sure never happens without the
  robber loosing.
  
  For ease of reference, we call the cop (or one of the cops) starting on $u$ the \textit{sentry,} the other the \textit{patrol.}
  Choose a vertex $v$ for $B$ as in Remark~\ref{rem:ex-coend}.  We now describe the winning strategy for the cops.  First of all, the patrol moves to $v$ on a
  finite path, while leaving the sentry on $u$.

  Now, cops play a winning strategy on the graph induced by $B$, chasing the projection of the robber on $G$ robber onto $B$, i.e., the vertex of $B$ which has
  smallest distance to the position of the robber in $G$.  By Lemma~\ref{lem:2-conn}, the cops either catch the robber or reach the cut-vertex $x$ of $G$ in $B$
  which is closest to the robber.  

  In the latter case, as soon as a cop reaches this vertex, he is put on sentry there and the other cop moves to $x$.  Letting $H'$ denote the connected component
  of $G\setminus\{x\}$ containing the robber and $H$ the subgraph of $G$ induced by $\{x\} \cup V(H')$, by induction, we can catch the robber in the graph $H$ with
  exit $x$.  This completes the proof of the theorem.
\end{proof}

We now give the proof of Lemma~\ref{lem:2-conn}.  It is easier to prove the following stronger version:

\begin{lemma}\label{lem:path-like}
  Let $G$ be  a path-like series parallel graph, and let $u$, $v$ be a pair of ends for $G$.  Then $G$ is $\{u,v\}$-exit 2-copwin.
\end{lemma}
For $B$ a block of a graph $G$ we call the set of all vertices of $B$ which are not cut-vertices of $G$ the \textit{interior} of $B$.
\begin{proof}
  We prove by induction on the number of edges of $G$, the case $\#E(G)=1$ being trivial.
  
  If $G$ has only one block, then, either the cops have won without fighting, or deleting $u$ and $v$ from $G$ leaves connected components, one of which, say $H'$
  contains the robber.  Then, the subgraph $H$ of $G$ induced by $\{u,v\}\cup V(H')$ is a path-like series-parallel graph.  We apply the induction hypothesis on
  $H$ with set of exits $u$, $v$.  Clearly, this proves the claim for $G$.

  Now suppose that $G$ has more than one block.  We need some terminology.  We call the cop starting on $x$ the $x$-cop, for $x=u,v$.  The \textit{home vertex} of
  the $x$-cop is $x$, the cops' \textit{home blocks} are defined similarly.  The $x$-cop's \textit{opposite vertex} is the cut-vertex of the home block, and
  denoted by $x'$.  The \textit{current vertex} of the $x$-cop at the end of move $t$ is denoted by $c_x(t)$.  Here, the movement of the two cops counts as one
  move, so in move $t=1$ it is the cops' turn, in move $t=2$ the robber's, and so on.
  
  The cops are now instructed to pursue the following strategy, either until one of them has reached his opposite vertex while the robber is not in this cop's home
  block, or, of course, the cops or the robber win the game.  We denote by $d(\cdot,\cdot)$ the shortest-path distance in the graph.

  \smallskip%
  \paragraph{Strategy of $x$-cop:}
  \begin{enumerate}
  \item Let $c$ be my current vertex, and $r$ the robber's current vertex.
  \item\label{enum:towards-home}%
    If $d(x,r) \le d(x,c)$, take one step on a shortest path from $c$ to $x$, which possibly requires to pass.
  \item\label{enum:towards-opp}%
    Otherwise take one step on a shortest path from $c$ to $x'$, which possibly requires to pass.
  \end{enumerate}

  We prove the following sequence of observations.
  \begin{description}
  \item[\it Claim 1.] The robber can never reach an exit unless being caught.  More precisely, after the cops move, the distance between the cop and his home
    vertex $x$ is at most the distance between the robber and $x$.
  \item[\it Claim 2.] A cop is always on a shortest path between his home vertex and his opposite vertex. 
  \item[\it Claim 3.] If the condition in (\ref{enum:towards-home}) is satisfied, then the robber is in the interior of the cop's home block.
  \item[\it Claim 4.] Define $\phi(t):=d(c_u(t),u')+d(c_v(t),v')$.  $\phi$ is decreasing (i.e., non-increasing).
  \item[\it Claim 5.] If $t$ is odd and $\phi(t)=\phi(t+2)$ then, after move $t$, the robber is in the interior of one cop's home block.
  \item[\it Claim 6.] At least one of the two robbers will reach his opposite vertex while the robber is not in his home block.
  \end{description}

  
  \subparagraph{\it Proof of Claim 1.}  Clear, because of the condition in (\ref{enum:towards-home}).

  \subparagraph{\it Proof of Claim 2.}  Clear by the movement of the cops.

  \subparagraph{\it Proof of Claim 3.}  Holds because every point which is not in the interior of the cop's home block is at least 1 farther away from his home
  vertex than any point on a shortest path between the home vertex and the in the opposite vertex, where, by Claim~2, the cop is.

  \subparagraph{\it Proof of Claim 4.}  By Claim~3, since the robber can be in at most one of the two home blocks at the same time.

  \subparagraph{\it Proof of Claim 5.}  The condition on $\phi$ implies that one of the functions $d(c_v(s),v')$, $d(c_u(s),u')$ increases from $s:=t$ to $s:=t+2$,
  while the other decreases.  The claim now follows from Claim~3.

  \subparagraph{\it Proof of Claim 6.}  By Claim~1, if this never happened, the game would continue infinitely.  Since the sequence $\phi$ is decreasing an
  non-negative, it must become stationary at some value.  This implies that the robber stays in one of the two home blocks for the rest of the game.  By Claim~5,
  this implies that the robber stays in the interior of one of the cops' home block for this time.  But then the other cop will reach his opposite vertex while the
  robber is not in his home block, a contradiction.

  \paragraph{\it Completion of the proof.}
  By Claim~6, we can assume that one of the two cops, say $u$ has cleared his home block from the robber, and it is the robbers turn.  By the second sentence in
  Claim~1, if we send the $v$-cop to his home vertex on a shortest path, he will reach it before the robber could.  We have thus reduced the number of blocks in
  the path-like series-parallel graph by one, and conclude by induction.
\end{proof}

\input{bibliography}
\end{document}

%% file: defs.amsart.tex
\theoremstyle{plain}
\newtheorem{theorem}{Theorem}

\newtheorem*{proposition*}{Proposition}

\newtheorem{lemma}[theorem]{Lemma}

\newtheorem*{theorem*}{Theorem}
\newtheorem*{lemma*}{Lemma}
\newtheorem*{conjecture*}{Conjecture}

\theoremstyle{definition}

\newtheorem*{exercise*}{Exercise}

\theoremstyle{remark}
\newtheorem{remark}[theorem]{Remark}
\newtheorem*{remark*}{Remark}
\newtheorem{remsTh}[theorem]{Remarks}

\newcommand{\subclass}[1]{}

\newcommand{\enumTi}[1]{\renewcommand{\theenumi}{#1}}

\input{defs.common}


%% file: defs.common.tex
\renewcommand{\em}{\sl}

\newlength{\algotabbingwidth}
